\documentclass{llncs}

\usepackage{upquote}
\usepackage[ngerman,main=english]{babel}
\usepackage{regexpatch}
\begin{document}

\title{PROOF OF  UNION - CLOSED SETS CONJECTURE}	
	\author{Vladimir Blinovsky\inst{1,2} \and Llohann D Speranca\inst{1}}
	\institute{
		Universidade Federal de São Paulo (UNIFESP). \\
		Campus Sao Jose dos Campos. Instituto de Ciencia e Tecnologia (ICT),
		Brazil\inst{1},\\
		Institute for Information Transmission Problems, \\
		B. Karetnyi 19, Moscow, Russia\inst{2}\\ 
\email{
		vblinovs@yandex.ru,\inst{1}\\
		lsperanca@gmail.com\inst{2}}}

	\maketitle

\begin{abstract}
We prove Union- Closed sets  conjecture
\end{abstract}
 \bigskip

The union- closed sets conjecture posed by Peter Frankl in 1979. There is an article in wikipedia on the URL~\cite{0}, devoted to this problem, see also~\cite{4} and List of unsolved mathematical problems~\cite{00}. A family of sets ${\cal A}\subset 2^{[n]}$ is said to be union - closed if the union of any two set from the family remains in the family. The conjecture states that for any union - closed family of subsets of finite set, other than family consisting only the empty set, there exists an element that belong to at least half of the sets from the family. In~\cite{5} Bruhn and Schaudt gave the nice survey of the result concerning the problem. 

The conjecture has been proven for many special cases. It is known to be true for families of at most 46 sets~\cite{1}, for $n\leq 11$~\cite{2}, for families of sets in which the smallest set has one or two elements~\cite{3}.

We use the natural bijection between $2^{[n]}$ and $\{ 0,1\}^n$ and don't make difference between these two sets. Considering natural embedding $\{ 0,1\}^n \to R^n$ we note, that arbitrary subset of $\{ 0,1\}^n$ can be defined by finite ($N$) number of inequalities
\begin{eqnarray}
\label{e11}
&&{\cal A}=\left\{ x\in \{ 0,1\}^n :\  \langle\omega_j,x\rangle > \delta_j (x ) ,\ j\in [N]\right\} ,
\end{eqnarray}
where $\sum_{q=1}^{n}\omega_{j,q}=C_j$  can be chosen as arbitrary (up to sign) given constants.

Define function 
$$
\varphi (\{\omega\}, x )\stackrel{\Delta}{=}\frac{1}{(2\pi )^{N/2}}\prod_{j=1}^N \int_{-\infty}^{(\langle\omega_j ,x \rangle-\delta_j (x) )/\sigma_j (m)}e^{-\xi^2 /2}d\xi\to\left\{\begin{array}{ll}
1,&  x\in{\cal A};\\
0,& \hbox{otherwise}
\end{array}
\right.
$$
as $\sigma_{j}(m)\to 0$ and $\sigma_j (m)\to 0$ as $m\to \infty$.

In the proof we use "Symmetrical Smoothing Method" which before was used in the proof of MMS conjecture by the first author~\cite{6}. 
Hence
$$
\biggl||{\cal A}|-\sum_{x\in 2^{[n]} \setminus\bar{0}}\frac{1}{(2\pi )^{N/2}}\prod_{i=1}^N \int_{-\infty}^{(\langle\omega_j ,x\rangle -\delta_j (x))/\sigma_j (m)}e^{-\xi^2 /2}d\xi\biggr|\to 0,$$
as $m\to \infty .$  Here we exclude all zero vector from summation over $x$. This is why $x\in 2^{[n]}\setminus\bar{0}$.

 Define 
 \begin{eqnarray}\label{d1} 
 &&R(\{\omega\}  ) = \sum_{x\in 2^{[n]}\setminus\bar{0}, x_{1}=0}\varphi (\{\omega\}, x)-\sum_{x\in 2^{[n]}\setminus\bar{0}, x_{1}=1}\varphi (\{\omega\}, x) ;\\
 \nonumber
 && M(\{\omega\} ) =\sum_{j=1}^N \sum_{x\in 2^{[n]}\setminus\bar{0}}\int_{-\infty}^{(\langle \omega_j  ,x \rangle -\delta_j (x))/\sigma_j (m)}e^{-\xi^2/2}d\xi \cdot \int_{-\infty}^{-(\langle  \omega_j ,x \rangle -\delta_j (x))/\sigma_j (m)}e^{-\xi^2/2}d\xi ; \\ \nonumber 
&& S_{\ell}(\{\omega\} ) = \sum_{x\in 2^{[n]}\setminus\bar{0},\ x_{\ell}= 1}\varphi (\{\omega\}, x )-\sum_{x\in 2^{[n]}\setminus\bar{0},\ x_{1}=1}\varphi (\{\omega\}, x),\ \ell =2,\ldots ,n; \\ \nonumber
&& L(\{\omega\} )=\sum_{x\neq x^\prime\in 2^{[n]}\setminus\bar{0}}\varphi (\{\omega\}, x )\cdot \varphi (\{\omega\}, x^\prime)\left(1-  \varphi \left(\{\omega\}, x\bigcup x^\prime \right)\right). \nonumber
\end{eqnarray}
W.l.o.g. we can also fix $\sum_{p=1}^n \omega_{j,p}=C_j =Const,\ j\in [N].$
\bigskip

Optimization Problem is to find
\begin{equation}
\label{eq1}
\max  R(\{\omega\} )
\end{equation}
when
\begin{eqnarray}\label{kk}
&& M(\{\omega\} )\leq  e^{-\Delta/\max_{i}\sigma_i};\  \Delta >0; \label{e0} \  S_{\ell}(\{\omega\} )\leq \frac{1}{2};\  L(\{\omega\} )\leq \frac{1}{2}\label{el2},
\end{eqnarray}
$m \to 0$. 

Value $M(\{\omega\})$ indicate the bound $|\langle \omega_j, x\rangle -\delta_j|>D$ for some universal constant $D=D(\Delta )>0$.
 Values $R(\{\omega\} )$ approximate the difference between the degree of vertex $0$ of ${\cal A}$ and degree of vertex $1$. Value $S_{\ell}(\{ \omega\} )$ approximate the difference between  degrees of vertex $\ell$ and  $1$. Value $L(\{\omega\} )$ is less that $\frac{1}{2}$ and this indicate that family ${\cal A}$ is union-closed
set.

 Fix $j\in [N]$ and define
 $$
 A(x,\{\omega\}, j)=\frac{1}{(2\pi)^{N/2}}\prod_{i=1,\ i\neq j }^N \int_{-\infty}^{(\langle x,\omega_i \rangle -\delta_i (x))/\sigma_i (m)} e^{-\frac{\xi^2}{2}}d\xi.
 $$
 We have 
 \begin{equation}
 \label{3w}
 \varphi (\{\omega\}, x )\stackrel{\Delta}{=}A(x,\{\omega\},j)\int_{-\infty}^{(\langle x,\omega_i \rangle -\delta_j (x))/\sigma_j (m)} e^{-\frac{\xi^2}{2}}d\xi .
 \end{equation}
 
 We have partial derivative 
  \begin{eqnarray} \label{c3}
&&(R(\{\omega\} ))^\prime_{\omega_{j,p}} =\frac{1}{\sigma_j}(A_1 -A_2), 
  \end{eqnarray}
  where
  \begin{eqnarray}
    && A_1 =\sum_{x\in 2^{[n]},\ x_{1} =0,\ x_p=1,\ x_{n} =0}A(x,\{\omega\},j)e^{-\frac{(\langle\omega_j ,x\rangle -\delta_j (x))^2}{2\sigma_j^2(m)} } \\ \nonumber
  &&+\sum_{x\in 2^{[n]},\ x_{1} =1,\ x_{p}=0,\ x_{n} =1} A(x,\{\omega\},j)e^{-\frac{(\langle \omega_j ,x\rangle -\delta_j (x))^2}{2\sigma_j^2 (m)} },\\ 
  &&A_2 =\sum_{x\in 2^{[n]},\ x_{1} =0,\ x_{p}=0,\ x_{n} =1}A(x,\{\omega\},j)e^{-\frac{(\langle \omega_j ,x\rangle -\delta_j (x))^2}{2\sigma_j^2 (m)} }\label{c4}\\
  &&+\sum_{x\in 2^{[n]},\ x_{1} =1,\ x_{p}=1,\ x_{n} =0}A(x,\{\omega\},j)e^{-\frac{(\langle \omega_j ,x\rangle -\delta_j (x))^2}{2\sigma_j^2 (m)} },\ j\in [N]. \nonumber
 \end{eqnarray}

 We find the unconditional extremum of the function from $R(\{\omega\})$ without taking into account conditions~(\ref{e0}) as the solution of the equations
\begin{equation}
\label{e110}
(R(\{\omega\} ))^\prime_{\omega_{j,p}} =0,\  \sum_{q=1}^n \omega_{j,q}=C_j ;\ p\in [2,n-1].
\end{equation}

Denote
 \begin{eqnarray*}&&
x(k) =\hbox{arg}\min_{x\in 2^{[n]}\setminus\left\{\bigcup_{i=1}^{k-1} x(i)\right\},\ \{ x_1 =1,\ x_p =0,\ x_n =1\},\ \{ x_1=0,\ x_p =1,\ x_n =0\}} (\langle \omega_j,x\rangle  -\delta_{j}(x))^2,\\
&& \bar{x}(k) =\hbox{arg}\min_{x\in 2^{[n]}\setminus\left\{ \bigcup_{i=1}^{k-1} \bar{x}(i)\right\},\ \{ x_1 =1,\ x_p =1,\ x_n =0\},\ \{ x_1=0,\ x_p =0,\ x_n =1\}}  (\langle \omega_j,x\rangle  -\delta_{j}(x))^2,\ k\in[2^{n-2}].
\end{eqnarray*}

We have
$$
A(x(k),\{\omega\} , j) = e^{\frac{\ln (A(x(k)\{\omega\}, j))\sigma^2_j (m)}{\sigma^2_j (m)}},\ A(\bar{x}(k),\{\omega\} , j) = e^{\frac{\ln (A(\bar{x}(k),\{\omega\}, j))\sigma^2_j (m)}{\sigma^2_j (m)}}
$$
and $\ln (A(x(k),\{\omega\}, j))\sigma^2_j (m),\  \ln (A(\bar{x}(k),\{\omega\}, j))\sigma^2_j (m)\to 0$ as $m\to \infty$. 

We can choose $\delta_j (x(k),\sigma_j (m)),\  \delta_j (\bar{x}(k),\sigma_j (m))$ s.t. 
\begin{eqnarray}
&& \label{3e}
A(x(k),\{\omega\} , j)e^{-\frac{(\langle \omega_j ,x(k)\rangle  -\delta_{j}(x(k)))^2}{2\sigma^2_j (m)}} =  e^{-\frac{ (\langle \omega_j ,x(k) \rangle  -\delta_{j}(x(k),\sigma_j(m)))^2}{2\sigma^2_j (m)}}\\ && 
A(\bar{x}(k),\{\omega\} , j)e^{-\frac{(\langle \omega_j ,\bar{x}(k)\rangle  -\delta_{j}(\bar{x}(k)))^2}{2\sigma^2_j (m)}} =  e^{-\frac{ (\langle \omega_j ,\bar{x}(k) \rangle  -\delta_{j}(\bar{x}(k),\sigma_j(m)))^2}{2\sigma^2_j (m)}}
\end{eqnarray}
and
\begin{equation}
\label{445}
e^{-\frac{ (\langle \omega_j ,x(k) \rangle  -\delta_{j}(x(k),\sigma_j(m)))^2}{2\sigma^2_j (m)}}=e^{-\frac{ (\langle \omega_j ,\bar{x}(k) \rangle  -\delta_{j}(\bar{x}(k),\sigma_j(m)))^2}{2\sigma^2_j (m)}}
\end{equation}
and
$$
\lim_{m\to\infty}(\langle \omega_j ,x\rangle  -\delta_{j}(x,\sigma_j(m)))^2  =( \langle \omega_j ,x \rangle  -\delta_{j}(x))^2.
$$
Denote
\begin{eqnarray*} &&
a(k)=\langle w_j,x(k) \rangle -\delta_{j}(x(k),\sigma_j(m)),\\
&& 
b(k)=\langle w_j,\bar{x}(k)\rangle  -\delta_{j}(\bar{x}(k),\sigma_j(m)).
\end{eqnarray*}

From equations~(\ref{445}) on the first step of iteration $k=1$ it follows 
 \begin{equation}
 \label{f5}
 (\langle w_j, x(1)\rangle  -\delta_{j} (x(1),\sigma_j(m)) )^2 = (\langle w_j,\bar{x}(1)\rangle -\delta_{j}(\bar{x}(1),\sigma_j(m)) )^2. 
 \end{equation}

Solution of this quadratic equation is the pair of equations
\begin{eqnarray}
\label{ww1}&&
\langle w_j,x (1)\rangle -\delta_{j}(x(1),\sigma_j(m)) - \langle w_j,\bar{x}(1)\rangle  -\delta_{j}(\bar{x}(1),\sigma_j(m)) =0,\\
\label{ww2}&& \langle w_j,x(1) \rangle -\delta_{j}(x(1),\sigma_j(m))  + \langle w_j,\bar{x}(1)\rangle -\delta_{j}(\bar{x}(1),\sigma_j(m))=0.
\end{eqnarray}

According to equalities~(\ref{ww1}),~(\ref{ww2}) we have two possibilities $a(1)=\pm b(1)$.   For proper choice of $a(1), b(1)$ we denote
$$
B(1)\in \left\{ \int_{-\infty}^{ a(1)/\sigma_j (m)} e^{-\xi^2 /2} d\xi \mp \int_{-\infty}^{\pm a(1)/\sigma_j (m)} e^{-\xi^2 /2} d\xi \right\} .
$$
We have $(B(1))^{\prime}_{\omega_{j,p}}=0$.

We establish the $i$-th step of iteration, which we describe as follows
$$
F(i)=R(\{\omega\}) -\sum_{k=1}^{i-1} B(k),\ i\in[2,2^{n-2}+1],
$$
where 
\begin{equation}
\label{dq1}
B(k)\in \left\{ \int_{-\infty}^{a(k)/\sigma_j (m)} e^{-\xi^2 /2} d\xi \mp \int_{-\infty}^{\pm a(k)/\sigma_j (m)}e^{-\xi^2 /2} d\xi\right\} ,\ (B(k))^\prime_{\omega_{j,p}} =0.
\end{equation}
And set 
$$
(F(i))^\prime_{\omega_{j,p}}=0.
$$

From equations~(\ref{445}) on $i$-th step it follows 
 \begin{equation}
 \label{f5}
 (\langle w_j, x(i)\rangle  -\delta_{j} (x(i),\sigma_j(m)) )^2 = (\langle w_j,\bar{x}(i)\rangle -\delta_{j}(\bar{x}(i),\sigma_j(m)) )^2. 
 \end{equation}

Solution of this quadratic equation is the pair of equations
\begin{eqnarray}
\label{w1}&&
\langle w_j,x (i)\rangle -\delta_{j}(x(i),\sigma_j(m)) - \langle w_j,\bar{x}(i)\rangle  -\delta_{j}(\bar{x}(i),\sigma_j(m)) =0,\\
\label{w2}&& \langle w_j,x(i) \rangle -\delta_{j}(x(i),\sigma_j(m))  + \langle w_j,\bar{x}(i)\rangle -\delta_{j}(\bar{x}(i),\sigma_j(m))=0.
\end{eqnarray}

According to equalities~(\ref{w1}),~(\ref{w2}) we have two possibilities $a(i)=\pm b(i)$  and 
$$
B(i)\in \left\{ \int_{-\infty}^{a(i)/\sigma_j (m)} e^{-\xi^2 /2} d\xi \mp \int_{-\infty}^{\pm a(i)/\sigma_j (m)} e^{-\xi^2 /2} d\xi \right\} .
$$
We have $(B(i))^{\prime}_{\omega_{j,p}}=0$.

Iteration of taking derivatives using the same arguments leads to the system of equalities
\begin{equation}
\label{ka1}
(\langle \omega_j ,x(k) \rangle -\delta_{j}(x(k),\sigma_j(m)))^2=(\langle \omega_j,\bar{x}(k)\rangle -\delta_{j}(\bar{x}(k),\sigma_j(m)))^2,\ k\in [2^{n-2}].
\end{equation}
 Summation of these equalities over $k\in [2^{n-2}]$ lead to the equalities
\begin{eqnarray*}
&&
\sum_{x\in 2^{[n]},\ \{ x_1 =1,\ x_p =0,\ x_n =1\},\ \{ x_1=0,\ x_p =1,\ x_n =0\}}(\langle \omega_j ,x\rangle -\delta_{j}(x,\sigma_j(m)) )^2 \\
&&=\sum_{x\in 2^{[n]},\ \{ x_1 =1,\ x_p =1,\ x_n =0\},\ \{ x_1=0,\ x_p =0,\ x_n =1\}}(\langle \omega_j ,x \rangle -\delta_{j}(x,\sigma_j(m)) )^2 ,\ p\in [2,n-1]
\end{eqnarray*}
or
\begin{eqnarray}\label{t5}&& \\
&& \nonumber
\sum_{x\in 2^{[n-3]}} \left(\omega_{j,p} +\sum_{q\neq 1,p,n}\omega_{j,q}x_q-\delta_{j}(x,\sigma_j(m))\right)^2    + \sum_{x\in 2^{[n-3]}} \left(\omega_{j,1} +\omega_{j,n}+\sum_{q\neq 1,p,n}\omega_{j,q}x_q-\delta_{j}(x,\sigma_j(m))\right)^2\\ \nonumber
&& 
=\sum_{x\in 2^{[n-3]}} \left(\omega_{j,n} +\sum_{q\neq 1,p,n}\omega_{j,q}x_q-\delta_{j}(x,\sigma_j(m))\right)^2  + \sum_{x\in 2^{[n-3]}} \left(\omega_{j,1} +\omega_{j,p}+\sum_{q\neq 1,p,n}\omega_{j,q}x_q-\delta_{j}(x,\sigma_j(m))\right)^2.
\end{eqnarray} 
Equality~(\ref{t5}) after routine transformations can be reduced to equality 
\begin{equation}
\label{ed1}
\omega_{j,p}=\omega_{j,n},\ p\in [2, n-1].  
\end{equation}

For $x\in \{0,1\}^n$ denote $|x|=\sum_{i=1}^n x_i .$
Now consider ${\cal A}_0 =\{ 1\not\in x \}\in{\cal A}\}$. Because $\omega_{j,p}=\omega_{j,n},\ p>1,\ j\in [N]$ and ${\cal A}$ is union-closed sets family, we have equality
$$
{\cal A}_0 =\{ x\in 2^{[2 ,n]}, |x|\geq \lambda\} 
$$
for some $\lambda$. Set ${\cal A}_1 ={\cal A}\setminus{\cal A}_0$ is nonempty, because $[n]\in{\cal A}_1$, and
$$
{\cal A}_1=\{ x\in 2^{[n]}:\ 1\in x,\ |x|\geq\lambda_1\}
$$
for some $\lambda_1$.
Because  for ${\cal A}_i  =\{ x\in{\cal A}:\ i\in x\}$ we have $|{\cal A}_i |\leq |{\cal A}_1 |$, it follows that ($i>1$)
\begin{eqnarray*}
|{\cal A}_i |&=&\sum_{j=\lambda}^{n-1}{n-2\choose j-1}+\sum_{j=\lambda_1}^{n}{n-2\choose j-2}\leq | {\cal A}_1 |\\ &=&\sum_{j=\lambda_1}^n {n-1\choose j-1}=\sum_{j=\lambda_1}^{n-1}{n-2\choose j-1}+\sum_{j=
\lambda_1}^{n}{n-2\choose j-2} .
\end{eqnarray*}
From here it follows that $\lambda_1 \leq \lambda $ and hence $2| {\cal A}_1| > |{\cal A}|$.

This completes  the proof of the conjecture.
\bigskip

{\bf Acknowledgment}
\bigskip

Research supported by the Sao Paulo Research Foundation (FAPESP), project no. 2012/13341-8 and  NUMEC/USP, project MaCLinC-USP 2013/07699-0 and visiting position in Unifesp, Extrato de contrato N114/2019, N20/2021, Processo N23089.101560/ 2018-91. First author would like to express gratitude to Unifesp and USP, where he started this work and specially to Prof.Y.Kohayakawa.

\end{document}